\documentclass{amsart}
\usepackage{amsmath,amsthm,amssymb,amscd,mathabx,accents,enumitem,mathtools,mathrsfs}
\usepackage{amsfonts}
\usepackage{rotating}
\usepackage{euscript}
\usepackage{pst-node}
\usepackage{epsfig,verbatim}
\usepackage{tikz-cd}
\def\be{\begin{equation}}
\def\ee{\end{equation}}

\def\C{{\mathbb C}} 
\def\f{\EuScript}

\def\P{{\mathbb P}}

\def\e{\eqref}
\def\phi{{\varphi}}

\def\v{{\varepsilon}} 
\def\tt{\widetilde}
\def\deg{{\rm deg\,}} 

\def\Aut{{\rm Aut}}

\def\mult{{\rm mult}}
 
\def\Ker{{\rm Ker\,}}

\def\mod{{\rm mod\ }}

\def\bp{\begin{proposition}}
\def\ep{\end{proposition}}

\def\bt{\begin{theorem}}
\def\et{\end{theorem}}
\def\br{\begin{remark}}
\def\er{\end{remark}}
\def\be{\begin{equation}}
\def\bee{\begin{equation*}}
\def\l{\label}

\def\e{\eqref}

\def\ee{\end{equation}}
\def\eee{\end{equation*}}
\def\bl{\begin{lemma}}
\def\el{\end{lemma}}
\def\bc{\begin{corollary}}
\def\ec{\end{corollary}}
\def\pr{\noindent{\it Proof. }}

\def\bd{\begin{definition}}
\def\ed{\end{definition}}

\mathtoolsset{showonlyrefs}

\newtheorem{theorem}{Theorem}[section]
\newtheorem{lemma}[theorem]{Lemma}
\newtheorem{definition}[theorem]{Definition}
\newtheorem{corollary}[theorem]{Corollary}
\newtheorem{proposition}[theorem]{Proposition}

\theoremstyle{definition}

\newtheorem{example}[theorem]
{Example}
\newtheorem{remark}[theorem]{Remark}

\begin{document}
\title{On symmetries of iterates of rational functions
}
\author{Fedor Pakovich}
\thanks{
This research was supported by ISF Grant No. 1092/22}
\address{Department of Mathematics, Ben Gurion University of the Negev, Israel}
\email{
pakovich@math.bgu.ac.il}

\begin{abstract} Let $A$  be  a rational function  of degree $n\geq 2$. Let us denote by 
$ G(A)$ the group of M\"obius transformations $\sigma$ such that $ A\circ \sigma=\nu_{\sigma} \circ A$ for some  M\"obius transformations $\nu_{\sigma}$, and by $\Sigma(A)$ and $\Aut(A)$ the subgroups of $ G(A)$ consisting of  $\sigma$ such that $ A\circ \sigma= A$ and  $ A\circ \sigma= \sigma \circ A$, correspondingly. In this paper, we study sequences of the above groups arising from iterating 
 $A$.  In particular, 
we show that if $A$ is not conjugate to $z^{\pm n},$ then   the orders of the groups   $ G(A^{\circ k})$, $k\geq 2,$ 
are finite and uniformly bounded in terms of  $n$ only.   
We also prove a number of results 
about the groups  $\Sigma_{\infty}(A)=\cup_{k=1}^{\infty} \Sigma(A^{\circ k})$ and $\Aut_{\infty}(A)=\cup_{k=1}^{\infty} \Aut(A^{\circ k})$, which are especially interesting from the dynamical perspective. 

\end{abstract}

\maketitle

\section{Introduction}
Let $A$ be a rational function of degree  $n\geq 2$.  
In this paper, we study a variety of different   subgroups of $\Aut(\C\P^1)$ related to $A$, and more generally to a dynamical system defined by iterating  $A$. Specifically, let us define $\Sigma(A)$ and  $\Aut(A)$ 
as the groups of M\"obius transformations $\sigma$ such that 
$ A\circ \sigma=A$ and 
$ A\circ \sigma= \sigma \circ A,$ correspondingly.
Notice that elements of $\Sigma(A)$ permute points of any fiber 
of $A$, and more generally of any fiber of $A^{\circ k},$ $k\geq 1,$ while  elements of $\Aut(A)$ permute 
fixed points of $A^{\circ k},$ $k\geq 1$. Since any M\"obius
transformation is defined by its values at any three points, this implies in particular that the groups $\Sigma(A)$ and $\Aut(A)$ are finite and therefore belong to the well-known list $A_4,$ $S_4,$ $A_5,$ $C_l$, $D_{2l}$  of finite subgroups of  $\Aut(\C\P^1)$.  

The both groups $\Sigma(A)$ and $\Aut(A)$ are subgroups of the group 
$ G(A)$ defined as the group of M\"obius transformations $\sigma$ such that \be \l{eblys} A\circ \sigma=\nu_{\sigma} \circ A\ee for some  M\"obius transformations $\nu_{\sigma}$.  It is easy to see that  $ G(A)$ is indeed a group, and that  $\nu_{\sigma}$ is defined in a unique way by $\sigma$.  Furthermore,  the map 
\be \l{homic}\gamma_A:\sigma \rightarrow \nu_{\sigma}\ee is a homomorphism from   $ G(A)$ to the group 
$\rm{Aut}(\C\P^1)$, whose kernel coincides with  $\Sigma(A)$. We will denote the image of $\gamma_A$ by $\widehat G(A).$ 
It was shown in the paper \cite{fin} that unless \be \l{for} A=\alpha\circ z^n\circ \beta\ee for some $\alpha,\beta\in \Aut(\C\P^1)$ the group $G(A)$ is also finite and its order is  bounded in terms of degree of $A$.

In this paper, we study the dynamical analogues 
of the groups $\Sigma(A)$ and $\Aut(A)$ defined by the formulas 
$$\Sigma_{\infty}(A)=\bigcup_{k=1}^{\infty} \Sigma(A^{\circ k}), \ \ \ \ \Aut_{\infty}(A)=\bigcup_{k=1}^{\infty} \Aut(A^{\circ k}).$$ 
Since 
\be \l{ospes} \Sigma( A)\subseteq \Sigma( A^{\circ   2})\subseteq \Sigma( A^{\circ   3})\subseteq \dots  \,\subseteq \Sigma( A^{\circ   k})\subseteq \dots \ ,\ee
and   \be \l{xre} \Aut(A^{\circ k})\subseteq  \Aut(A^{\circ r}),\ \ \  \  \ \Aut(A^{\circ l})\subseteq  \Aut(A^{\circ r})\ee 
for any common multiple $r$ of $k$ and $l$,   
the sets $\Sigma_{\infty}(A)$ and $\Aut_{\infty}(A)$ are  groups. 
 While it is not clear a priori that the groups  $\Sigma_{\infty}(A)$ and $\Aut_{\infty}(A)$  
are finite, for $A$ not conjugated to $z^{\pm n}$ their finiteness can be deduced from the theorem of Levin (\cite{lev}, \cite{lev2}) about rational functions sharing the measure of maximal entropy. 
However, the Levin theorem  does not permit to describe the groups $\Sigma_{\infty}(A)$ and $\Aut_{\infty}(A)$  or to estimate their orders, and the main goal of this paper is  to prove some results in this direction. 
More generally,   we study the totality of the groups  $ G(A^{\circ k})$, $k\geq 1,$
defined by iterating $A.$ 

Our main result about the groups  $ G(A^{\circ k})$, $k\geq 1,$ can be formulated as follows.

\bt \l{sig0} Let $A$ be a rational function of degree $n\geq 2$ that is not conjugate to $z^{\pm n}.$
Then   the orders of the groups   $ G(A^{\circ k})$, $k\geq 2,$ 
are finite and uniformly bounded in terms of  $n$ only. 
\et

In addition to Theorem \ref{sig0}, we prove a number of more precise results about the groups $\Sigma_{\infty}(A)$ and $\Aut_{\infty}(A)$ allowing us in certain cases to calculate these groups  explicitly. For a rational function $A$, let us  denote by $c(A)$ the set of its critical values. Our main result concerning the  groups $\Aut_{\infty}(A)$ is following.

\bt \l{main} Let $A$ be a rational function of degree $n\geq 2$ that is not conjugate to $z^{\pm n}.$ 
Then  the group $\Aut_{\infty}(A)$ is finite and  its order is bounded in terms of $n$ only.  Moreover,  every $\nu\in \Aut_{\infty}(A)$ maps the set 
$c(A)$ to the set $c(A^{\circ 2})$.
\et

Notice that since  M\"obius transformations $\nu$ such that  \be \l{st} \nu\big(c(A)\big)\subseteq c(A^{\circ 2})\ee can be described explicitly, Theorem \ref{main} provides us with a concrete subset of $\Aut(\C\P^1)$ containing  the group  $\Aut_{\infty}(A)$.

To formulate our main results concerning groups $\Sigma(A)$, let us introduce some definitions. 
Let $A$ be a rational function. Then a rational function $\widetilde A$ is called an {\it elementary transformation} of $A$
if there exist rational functions $U$ and $V$ such that   \be \l{ey} A=U\circ V\ \ \ \ {\rm and} \ \ \ \ \widetilde A=V\circ U.\ee We say that rational functions 
 $A$ and $A'$ are  {\it equivalent} and write $A\sim A'$  if there exists 
a chain of elementary transformations between $A$ and $A'$.
Since for any M\"obius transformation $\mu$ the equality 
\be \l{lus} A=(A\circ \mu^{-1})\circ \mu\ee holds, 
the equivalence class $[A]$ of a rational function $A$ is a union of conjugacy classes.
Moreover, by the results of the papers \cite{rec}, \cite{fin}, the number of conjugacy classes in $[A]$ is finite, unless $A$ is a flexible Latt\`es map.

In this notation, our main result about the groups $\Sigma_{\infty}(A)$ is following. 

\bt \l{sug}
Let $A$ be a rational function of degree $n\geq 2$ that is not conjugate to $z^{\pm n}.$ Then the order of the group  $\Sigma_{\infty}(A)$ is finite and  bounded in terms of $n$ only. Moreover, 
for every $\sigma\in\Sigma_{\infty}(A)$ the relation $A\circ \sigma\sim A$ holds. 
\et

Notice that  
in some cases Theorem \ref{sug} permits to describe the group $\Sigma_{\infty}(A)$  completely. 
Specifically, assume that   $A$ is {\it indecomposable}, that is, cannot be represented as a composition of two rational functions of degree at least two. In this case, 
 the number of conjugacy classes in the equivalence class $[A]$ obviously is equal to one, and Theorem \ref{sug} yields the following statement. 

\bt \l{sig4} Let $A$ be an indecomposable rational function of degree $n\geq 2$ that is not conjugate to $z^{\pm n}.$  Then  $\Sigma_{\infty}(A)=\Sigma(A)$, 
whenever the group $\widehat G(A)$ is trivial.  
 Moreover,  the group $\Sigma_{\infty}(A)$ is  trivial, whenever $G(A)=\Aut(A)$. 
\et

Notice that Theorem \ref{sig4} implies in particular that if $A$ is indecomposable and the group $G(A)$ is trivial, then 
$\Sigma_{\infty}(A)$ is also trivial.

Finally, along with the groups  $G(A^{\circ k}),$ $k\geq 1,$  we consider their ``local'' versions.
Specifically, let $z_0\in \C\P^1$  be a  fixed point of  $A$. For  a point  $z_1\in \C\P^1$ distinct from $z_0$, we define $G(A,z_0,z_1)$ as the subgroup of $G(A)$  consisting of M\"obius transformations $ \sigma$ such that  $ \sigma(z_0)=z_0$ and $ \sigma(z_1)=z_1$. 
For these groups, we prove the following statement.

\bt \l{seq0}  Let $A$ be a rational function of degree  $n\geq 2$ that is not conjugate to $z^{\pm n}.$
 Assume that $z_0\in \C\P^1$ is  a  fixed point of $A$,  and   $z_1\in \C\P^1$ is a point distinct from $z_0$. 
 Then $G(A^{\circ k},z_0,z_1)$, $k\geq 1,$ are finite cyclic groups equal to each other.  
\et

Notice that every element $\sigma \in \Aut(A^{\circ k})$, $k\geq 1$,  
belongs to $G(A^{\circ 2k},z_0,z_1)$ for some $z_0,$ $z_1$. Indeed, 
 the equality 
\be \l{au} A^{\circ k}\circ \sigma=  \sigma \circ A^{\circ k}, \ \ \ \ k\geq 1, \ee  
implies that $A^{\circ k}$ sends the set of fixed points of $\sigma$ to itself. 
Therefore, at least one of these  points $z_0,$ $z_1$ is a fixed point of $A^{\circ 2k},$ and if $z_0$ is such a point, then  
 $\sigma \in  G(A^{\circ 2k},z_0,z_1).$
In view of this relation between $\Aut(A^{\circ k})$ and $G(A^{\circ 2k},z_0,z_1),$ Theorem \ref{seq0} allows us  in some cases to estimate the order of the group $\Aut_{\infty}(A)$ and even to describe this group explicitly.

The paper is organized as follows. In the second section,
we establish basic properties of the group $G(A)$  
and provide a method for its calculation. In the third section,   we briefly discuss relations between the groups $\Sigma_{\infty}(A)$, $\Aut_{\infty}(A)$ and the  measure of maximal entropy for $A$. In particular, we deduce the  finiteness  of these  groups from the results of Levin (\cite{lev}, \cite{lev2}).

In the fourth  section,  we prove Theorem \ref{main}. Moreover, we prove that  \eqref{st} holds for any M\"obius transformation $\nu$ that belongs to $\widehat  G(A^{\circ k})$ for some $k\geq 1.$ 
 In the  fifth section, using results about semiconjugate rational functions from the papers \cite{semi},  \cite{fin}, we prove Theorem \ref{sug} and Theorem \ref{sig4}. We also prove a slightly more general version of Theorem \ref{sig0}.  
Finally, in the sixth section, we deduce Theorem \ref{seq0} from the result of Reznick (\cite{rez}) about iterates of formal power series, and  provide some applications of Theorem \ref{seq0} concerning the groups  $\Aut_{\infty}(A)$ and $\Sigma_{\infty}(A)$.

\section{Groups $G(A)$}
Let $A$ be a rational function of degree $n\geq 2$, and $G(A),$ $\widehat G(A),$ $\Sigma(A),$ $\Aut(A)$ the groups defined in the introduction. 
Notice that if rational functions $A$ and $A'$ are related by the equality 
$$\alpha\circ  A\circ \beta=A'$$ for some $\alpha,\beta\in \Aut(\C\P^1)$, 
then 
\be \l{siom} G(A')= \beta^{-1}\circ G(A) \circ \beta,  \ \ \ \ \ \ \ \widehat G(A')=\alpha \circ \widehat G(A) \circ \alpha^{-1}.\ee 
In particular, the groups $G(A)$ and $G(A')$ are isomorphic.
 Notice also that since
\be \l{isom} \widehat G(A)\cong G(A)/\Sigma(A),\ee the equality 
\be \l{burn} \vert  G(A)\vert =\vert \widehat G(A)\vert\vert \Sigma(A)\vert\ee holds 
whenever  the groups involved are finite.

\bl \l{meb} Let  $A$ be a rational function of degree $n\geq 2$. Then the following 
statements are true.

\begin{enumerate}[label=\roman*)]
\item For every  $z\in \C\P^1$ and  $\sigma\in G(A)$ the multiplicity of $A$ at $z$ is equal to the multiplicity of $A$ at $\sigma(z).$
\item For every  $c\in \C\P^1$ and  $\sigma\in G(A)$ the fiber $A^{-1}\{c\}$ is mapped by $\sigma$ to the fiber $A^{-1}(\nu_{\sigma}(c))$. 
\item Every 
$\nu\in \widehat G(A)$ maps $c(A)$ to $c(A).$
\end{enumerate}

\el 
\pr Since $\eqref{eblys}$ implies that 
$$\mult_{\sigma(z)}A\cdot\mult_z\sigma=\mult_{A(z)}\nu_{\sigma}\cdot \mult_{z}A $$
 the first statement follows from the fact that $\sigma$ and $\nu_{\sigma}$  are one-to-one. 

Further, it is clear that \eqref{eblys} implies 
$$\sigma^{-1}(A^{-1}\{c\})=A^{-1}(\nu_{\sigma}^{-1}\{c\}).$$
Changing now $\sigma^{-1}$ to $\sigma$ and taking into account that 
$\nu_{\sigma}^{-1}=\nu_{\sigma^{-1}}$, we obtain the second statement. 

Finally, the third statement follows from the second one, taking into account that 
$$\vert A^{-1}\{c\}\vert= \vert A^{-1}\{\nu_{\sigma}(c)\}\vert $$ since 
$\sigma$ is one-to-one, and that 
$c$ is a critical value of $A$ if and only $\vert A^{-1}\{c\}\vert <n$.  \qed 

We say that a rational function $A$ of degree $n\geq 2$ is {\it a quasi-power} if 
there exist $\alpha,\beta\in \Aut(\C\P^1)$ such that $$A= \alpha\circ z^n\circ \beta.$$
It is easy to see using Lemma \ref{meb} that  the group $G(z^n)$ consists of the transformations $z\rightarrow c z^{\pm 1}$, $c\in \C\setminus\{0\}$. Therefore, by \eqref{siom}, for any quasi-power $A$ 
the groups $G(A)$ and $\widehat G(A)$ are infinite. 

\bl \l{qua}  A rational function $A$ of degree $n\geq 2$  is a quasi-power if and only if it has only two critical values. If $A$ is a quasi-power, then $A^{\circ 2}$ is a quasi-power if and only if $A$ is conjugate to  $z^{\pm n}.$ 
\el 
\pr The first part of the lemma is well-known and follows easily from the Riemann-Hurwitz formula. To prove the second, we observe that the chain rule implies that the function 
$$A^{\circ 2}=\alpha\circ z^n\circ \beta\circ \alpha\circ z^n\circ \beta$$ has only two critical values if and only if  $\beta\circ \alpha$ maps the set $\{0,\infty\}$ to itself. Therefore, $A^{\circ 2}$ is a quasi-power if and only if $\beta\circ \alpha=cz^{\pm 1},$ $c\in \C\setminus \{0\},$ that is, if and only if 
$$A=\alpha\circ z^n\circ \beta=\alpha\circ z^n\circ cz^{\pm 1} \circ \alpha^{-1}=\alpha\circ c^nz^{\pm n}\circ \alpha^{-1}.$$  
Finally, it is clear that the last condition is equivalent  to the condition that $A$ is conjugate to $z^{\pm n}$.  \qed 

Let $G$ be a finite subgroup of $\Aut(\C\P^1)$. We recall that a rational function $\theta_G$ is called an {\it invariant function}  for $G$ if  the equality $\theta_G(x)=\theta_G(y)$  holds for $x,y\in \C\P^1$ if and only if there exists $\sigma\in G$ such that $\sigma(x)=y.$ Such a function always exists and is defined in a unique way up to the transformation $\theta_G\rightarrow \mu \circ \theta_G,$ where $\mu\in \Aut(\C\P^1)$. Obviously, $\theta_G$ has degree equal to the order of $G$. Invariant functions for finite subgroups of $\Aut(\C\P^1)$ were first found by Klein in his book \cite{klein}.

\bt \l{aut} 
Let $A$ be a rational function of degree $n\geq 2$. 
  Then  $\Sigma(A) $ is a finite group and $\vert  \Sigma(A)\vert $ is a divisor of $n$. Moreover, $\vert  \Sigma(A)\vert =n$ if and only if 
$A$ is an invariant function for $\Sigma(A).$ 
\et 
\pr Since for a finite subgroup $G$ of $\Aut(\C\P^1)$ the set of rational  functions $F$ such that $F\circ \sigma=F$ for every $\sigma \in G$ is a subfield of $\C(z),$ it follows easily from   
 the L\"uroth theorem that any such a function  $F$ is a rational function in $\theta_G$. Thus, $\deg F$ is divisible by $\deg\theta_G=\vert G \vert.$ In particular, setting $G=\Sigma(A)$, we see that the degree of $A$  is divisible by   $\vert  \Sigma(A)\vert$, and  $\deg A=\vert  \Sigma(A)\vert$ if and only if 
$A$ is an invariant function for $\Sigma(A).$  
\qed

The existence of invariant functions implies that for every finite subgroup $G$ of $\Aut(\C\P^1)$ there exist rational functions for which $\Sigma(A)=G$.  
Similarly, for every finite subgroup $G$ of $\Aut(\C\P^1)$
there exist  rational functions for which $\Aut(A)=G$. A description of such functions in terms of homogenous invariant polynomials for $G$ was obtained by Doyle and McMullen in \cite{dm}.  
Notice that rational functions with non-trivial automorphism groups are closely related to {\it generalized Latt\`es maps} (see \cite{lattes} for more detail).

The following result was proved in \cite{fin}. For the reader convenience we provide a  simpler   proof.

\bt \l{prim} Let $A$ be a rational function of degree $n\geq 2$ that is not a quasi-power. 
  Then the group $ G(A)$ is isomorphic to one of the five finite rotation groups of the sphere $A_4,$ $S_4,$ $A_5,$ $C_l$, $D_{2l}$, 
and	the order  of any element of $ G(A)$ does not exceed  $n.$ In particular, $\vert  G(A)\vert \leq \max\{60,2n\}.$
\et 

\pr 
Any  element of the group $\Aut(\C\P^1)\cong \rm{PSL}_2(\C)$ is conjugate either to $z\rightarrow z+1$ or to $z\rightarrow \lambda z$ for some $\lambda \in \C\setminus\{0\}.$ Thus, making the change 
\be \l{zve} A\rightarrow \mu_1\circ A\circ \mu_2, \ \ \ \sigma\rightarrow \mu_2^{-1}\circ \sigma\circ \mu_2,\ \ \ \nu_{\sigma}\rightarrow \mu_1\circ \nu_{\sigma}\circ \mu_1^{-1}\ee for convenient $\mu_1,$ $\mu_2\in\rm{Aut}(\C\P^1)$, without loss of generality we may assume that $\sigma$ and $\nu_{\sigma}$ in \eqref{eblys} have one of the two forms above.  

We observe first that 
the equality \be \l{j1} A(z+1)=\lambda A(z), \ \ \ \ \ \lambda\in \C\setminus\{0\},\ee 
is impossible. Indeed, if $A$ has a finite pole, then \eqref{j1} implies that $A$ has infinitely many poles. 
On the other hand, if $A$ does not have finite poles, then 
$A$ has a finite zero, and \eqref{j1} implies that $A$ has infinitely many zeroes. 
Similarly, the equality 
 \be \l{j2}A(z+1)=A(z)+1\ee is impossible
if $A$ has a finite pole. 
On the other hand, if $A$ is a polynomial of degree $n\geq 2$, then we obtain a contradiction comparing  the coefficients of
 $z^{n-1}$ on the left and the right sides of  equality \eqref{j2}.

 For the argument below, instead of considering $A$ as a ratio of two polynomials,  
it is more convenient to assume that $A$ is represented by its convergent 
 Laurent series at zero or infinity. 
Comparing for such a representation the free terms on the left and the right sides of the equality \be \l{egga} A(\lambda z)=A(z)+1, \ \ \ \lambda\in \C\setminus\{0\},\ee 
we conclude that this equality is impossible either. 
Thus, equality \eqref{eblys} for a non-identity $\sigma$  reduces to the equality 
\be \l{zore} A(\lambda_1z)=\lambda_2A(z),\ \ \ \lambda_1\in \C\setminus\{0,1\},\ \ \  \lambda_2\in \C\setminus\{0\}.\ee  
Comparing now coefficients on the left and the right sides of \eqref{zore} and taking into account that
 $A\neq  az^{\pm n}$, $a\in \C,$ by the assumption, 
we conclude that $\lambda_1$ is a root of unity.
Furthermore, if $d$ is the order of $\lambda_1$, then 
$\lambda_2=\lambda_1^r$ for some $0\leq r \leq d-1$, implying that  $A/z^r$ is a rational function in $z^d$. 
On the other hand, it is easy to see that if $ A=z^rR(z^d),$ where $R\in \C(z)$ and $0\leq r \leq d-1,$ 
then $d\leq  n,$ unless either $R\in \C\setminus\{0\}$ or $R=a/z$ for some $a\in \C\setminus\{0\}$. 
Since for such $R$ the function $A$ is a quasi-power, we conclude that  
 the order of $\lambda_1$ and hence the order of any element of  $ G(A)$ does not exceed $n$.

To finish the proof we only must show that  $G(A)$ is finite. By Lemma \ref{qua},  $A$ has at least three critical values. On the other hand, by 
Lemma  \ref{meb}, iii), every 
$\nu\in \widehat G(A)$ maps $c(A)$ to $c(A).$ Since any M\"obius transformation is defined by its values at any three points, this implies that $\widehat G(A)$ is finite. Since $\Sigma(A)$ is finite by Theorem \ref{aut},  
this implies that  $G(A)$ is finite because of the isomorphism \eqref{isom}.
 \qed 

\begin{remark} 
Using some non-trivial group-theoretic results about subgroups of $\rm{GL}_k(\C)$, one can deduce the finiteness of $G(A)$ directly from the fact that 
the order  of any element of $ G(A)$ does not exceed  $n.$ 
 Namely, the proof given in the paper \cite{fin} uses the Schur theorem (see e.g. \cite{cur}, (36.2)),  which states that any  finitely generated periodic subgroup of  $\rm{GL}_k(\C)$ has finite order. Alternatively, one can use the Burnside theorem (see e.g. \cite{cur}, (36.1)),  which states that any subgroup of $\rm{GL}_k(\C)$ of bounded period is finite. 
Indeed, assume that $G(A)$ is infinite. Then its lifting $\overline{G(A)}\subset \rm{SL}_2(\C)\subset \rm{GL}_2(\C)$ is  also 
infinite. On the other hand, if the order of any element of $ G(A)$ is  bounded by $N$, then 
 the order of any element of $\overline{G(A)}$ is  bounded by $2N$. The contradiction obtained proves the  finiteness of $G(A)$. 
\end{remark}

\bc \l{yc2} Let $A$ be a rational function of degree $n\geq 2$. Then 
$\Sigma(A)$ and $\Aut(A)$ are finite groups whose order does not exceed $\max\{60,2n\}.$
\ec 
\pr  If $A$ is a not a quasi-power, then the corollary follows from Theorem \ref{prim}. On the other hand, it is easy to see that if $A$ is a quasi-power, then 
the corresponding groups are cyclic groups of order $n$ and $n-1$ correspondingly. \qed

Let us mention the following specification of Theorem \ref{prim}. 

\bt \l{cyc} Let $A$ be a rational function of degree $n\geq 2$.  Assume that there exists a point $z_0\in \C\P^1$ 
such that the multiplicity of $A$ at $z_0$ is distinct from the multiplicity of $A$ at any other point $z\in \C\P^1$. Then $G(A)$ is a finite cyclic group, and  $z_0$  is a 
fixed point of its generator.
\et 
\pr 
It follows from the assumption that $A$ is not a quasi-power. Therefore, $G(A)$ is finite. 
Moreover, every element of $G(A)$ fixes $z_0$ by Lemma \ref{meb}, i). On the other hand, a unique finite subgroup of $\Aut(\C\P^1)$ whose elements share a fixed point is cyclic. \qed 

In turn, Theorem \ref{cyc} implies the following well-known corollary.

\bc \l{cyc2} Let $P$ be a polynomial of degree $n\geq 2$ that is not a quasi-power. Then $G(P)$ is a finite cyclic group generated by a polynomial.  
\ec 
\pr  Since $P$ is a not a quasi-power, the multiplicity of $P$ at infinity is distinct from 
the multiplicity of $P$ at any other point of $\C\P^1.$ Moreover, since every element of $G(P)$ fixes 
infinity, $G(P)$ consist of polynomials.
\qed

Notice  that functions $A$ of degree $n$ with $\vert G(A)\vert =2n$ do exist.   
Indeed, it is easy to see that  for any function of the from $$A=\frac{z^n-a}{az^n-1},  \ \ \ a\in \C\setminus\{0\},$$ the group $G(A)$ contains the dihedral group $D_{2n}$, generated by 
$$z\rightarrow \frac{1}{z}, \ \ \ \ \ z\rightarrow \v_nz,$$ where $\v_n=e^{\frac{2\pi i}{n}}.$ 
Thus,  for $n$ big enough, $G(A)=D_{2n}$, by Theorem \ref{prim}. 
On the other hand, for small $n$, functions $A$ of degree $n$ with $\vert G(A)\vert >2n$  
do exist as well (see for instance Example \ref{dadef} below).   

Lemma \ref{meb} provides us with a  method for practical calculation of $G(A)$, at least if the degree of $A$ is small enough. We illustrate it with the following example.

\begin{example}\l{exa1}
Let us consider  the function 
\be \l{seen} A=\frac{1}{8}\,{\frac {{z}^{4}+8\,{z}^{3}+8\,z-8}{z-1}}.\ee
One can check that $A$ has three critical values $1$, $9$, and $\infty$, and that 
$$ A-1= \frac{1}{8}\,{\frac {{z}^{3} \left( z+8 \right) }{z-1}}, \ \ \ \ \ \ \  A-9= \frac{1}{8}\,{\frac { \left( {z}^{2}+4\,z-8 \right) ^{2}}{z-1}}.$$
Since the multiplicities of $A$ at the preimages of   $1$, $9$, and $\infty$ are 
$$\mult_{0}A=3, \ \ \ \mult_{-8}A=1, \ \ \ \mult_{-2+2\sqrt{3}}A=2, \ \ \ \mult_{-2- 2\sqrt{3}}A=2,$$ and $$ \mult_{\infty}A=3, \ \ \ \mult_{1}A=1, $$
Lemma \ref{meb} implies that for any $\sigma \in G(A)$ 
 either 
\be \l{us1} \sigma(0)=0, \ \ \sigma(\infty)=\infty, \ \  \sigma(-8)=-8, \ \ \  \sigma(1)=1,\ee or 
\be \l{us2} \sigma(0)=\infty, \ \ \sigma(\infty)=0, \ \  \sigma(-8)=1, \ \ \  \sigma(1)=-8.\ee
Moreover, in addition, either 
\be \l{us3}  \sigma(-2+ 2\sqrt{3})=-2- 2\sqrt{3}, \ \ \  \sigma(-2- 2\sqrt{3})=-2+ 2\sqrt{3},\ee 
or 
\be \l{us4}  \sigma(-2+ 2\sqrt{3})=-2+ 2\sqrt{3}, \ \ \  \sigma(-2- 2\sqrt{3})=-2- 2\sqrt{3}.\ee

 Clearly, condition \eqref{us1} implies that $\sigma=z$, while the unique transformation satisfying \eqref{us2} is \be \l{fo1} \sigma=-8/z,\ee and this transformation satisfies \eqref{us3}.  
Furthermore, the corresponding $\nu_{\sigma}$ must satisfy 
$$\nu_{\sigma}(1)=\infty, \ \ \ \ \nu_{\sigma}(\infty)=1, \ \ \ \ \nu_{\sigma}(9)=9,$$ implying that 
\be \l{fo2} \nu_{\sigma}=\frac{z+63}{z-1}.\ee Therefore, \eqref{eblys} can hold only for $\sigma$ and $\nu_{\sigma}$ given by formulas \eqref{fo1} and \eqref{fo2}, and a direct calculation shows that \eqref{eblys} is indeed satisfied. Thus, the group $G(A)$ is a cyclic group of order two.

\end{example}

Notice that to verify whether a given M\"obius transformation  $\sigma$ belongs to $G(A)$ one can use the  Schwarz derivative. Let us recall that 
for a function $f$ meromorphic on a domain $D\subset\C$ the Schwarz derivative is defined by
\be \l{dsch}
S(f)(z)=\frac{f'''}{f'}-\frac{3}{2}\left(\frac{f''}{f'}\right)^2.
\ee
The characteristic property of the Schwarz derivative is that for two functions
 $f$ and $g$ meromorphic on $D$ the equality $S(f)(z)=S(g)(z)$ holds if and only if 
$g=\nu \circ f$ for some M\"obius transformation $\nu.$ 
%$$g\,=\,\frac{a\,f+b}{c\,f+d} \mbox{ for some } {\binom{a\,\ b}{c\,\ d}}\in {\mbox{GL}_2(\C)}.
%$$
Thus,  a  M\"obius transformation $\sigma$ belongs to $G(A)$ if and only if $$S(A)(z)=S(A\circ\sigma) (z).$$

We finish this section by another example of calculation of $G(A).$

\begin{example} \l{dadef} 
Let us consider the function 
$$B=-\frac{2z^2}{z^4+1}=-\frac{2}{z^2+\frac{1}{z^2}}.$$ It is easy to see that $\Sigma(B)$ contains the transformations $z\rightarrow -z$ and $z\rightarrow 1/z$, which generate  
the Klein four-group $V_4=D_4$, implying that $\Sigma(B)=D_4$ by Theorem \ref{aut}. Furthermore, it is clear that $G(B)$ contains the transformation $z\rightarrow iz$, implying that   
 $G(B)$ contains $D_8.$

The groups $A_4$, $A_5$,  and $C_l$ do not contain $D_8$. Therefore, if  $D_8$ is a proper subgroup of $G(B)$, then either  $G(B)=S_4$, or  $G(B)$ is a dihedral group containing an element $\sigma$ of order $k>4$, whose fixed points coincide with fixed points of $z\rightarrow iz$. The second case is impossible, since any M\"obius transformation $\sigma$ fixing $0$ and $\infty$  has the form  $cz$, $c\in \C\setminus\{0\},$ and it is easy to see that such $\sigma$ belongs to $G(B)$ if and only if it is a power of $z\rightarrow iz.$   On the other hand, a  direct calculation shows that for the  transformation $\mu=\frac{z+i}{z-i}$, generating together with $z\rightarrow iz$ and $z\rightarrow 1/z$ the group $S_4$,  equality \eqref{eblys}  holds for $\nu={\frac {-z+1}{-3\,z-1}}.$ 
Thus, $G(B)\cong S_4.$

\end{example}

\section{ Groups $\Sigma_{\infty}(A)$, $\Aut_{\infty}(A)$ and the  measure of maximal entropy}  
Let us recall that by the results of Freire, Lopes, Ma\~n\'e (\cite{flm}) and Lyubich (\cite{l}), 
 for every rational function $A$ of degree $n\geq 2$ there exists a unique probability measure $\mu_A$ on $\C\P^1$, which is invariant under $A$, has support equal to the Julia set $J_A$, and achieves maximal entropy 
$\log n$ among all $A$-invariant probability measures.

The measure $\mu_A$ can be described as follows. For $a\in \C\P^1$ 
let  $z_i^k(a),$ $i=1, \dots, n^k,$ be the roots of the equation
$A^{\circ k}(z) = a$ counted with multiplicity, and $\mu_{A,k}(a)$  
the measure defined by 
\be \l{mera} \mu_{A,k}(a)=\frac{1}{n^k}\sum_{i=1}^{n^k}\delta_{z_i^k(a)}.\ee Then 
for every  $a\in \C\P^1$ with two possible exceptions, the sequence  $\mu_{A,k}(a)$, $k\geq 1$,  converges in the weak topology to $\mu_A.$ 
%It follows from the characterization of $\mu_F$ as a limit of \eqref{mera}  that 
%This implies i
Notice that this description of $\mu_A$ implies that  $\mu_A=\mu_B$ whenever $A$ and $B$ share an iterate. 

The measure $\mu_A$
is characterized by the balancedness property that $$\mu_A(A(S))=\mu_A(S)\deg A\,$$ 
for any Borel set $S$ on which $A$ is injective.  Notice that for rational functions $A$ and $B$ the property to have the same measure of maximal entropy can be 
expressed also in algebraic terms (see \cite{lp}), leading to characterizations of such functions in terms of functional equations (see \cite{lp}, \cite{entr}, \cite{ye}).

The relations between the groups $\Sigma_{\infty}(A)$, $\Aut_{\infty}(A)$ and the  measure of maximal entropy are described by the following two statements.  

\bl \l{fst} Let $A$ be a rational function of degree $n\geq 2$. Then $\sigma\in \Aut_{\infty}(A)$ if and only if $A$ and $\sigma^{-1}\circ A\circ \sigma$
 have a common iterate. In particular, if $\sigma\in \Aut_{\infty}(A)$, then $A$ and $\sigma^{-1}\circ A\circ \sigma$
 share the  measure of maximal entropy.
\el
\pr The proof is trivial, given that rational functions sharing an iterate share a measure of maximal entropy.  \qed 

\bl \l{sst} Let $A$ be a rational function of degree $n\geq 2$. Then for every \linebreak $\sigma\in \Sigma_{\infty}(A)$ the functions $A$ and $A\circ \sigma$ share the measure of maximal entropy. 
\el
\pr The equality $$A^{\circ l}=A^{\circ l}\circ \sigma,\ \ \ l\geq 1,$$  implies that for any $k\geq l$ and $a\in \C\P^1$ the transformation $\sigma$ maps 
 the set of roots of the equation 
$A^{\circ k}(z) = a$ to itself. Thus, 
 for any set $S\subset \C\P^1$ we have 
$$\vert S\cap A^{-k}(a)\vert =\vert \sigma(S)\cap A^{-k}(a)\vert, \ \ \ k\geq l, \ \ \ a\in \C\P^1,$$ implying 
that  any $\sigma\in \Sigma_{\infty}(A)$ is $\mu_A$-invariant since $\mu_A$ is a limit of \eqref{mera}.

 Let now $S$ be a Borel set on which $A\circ \sigma$ is injective. Then $A$ is injective on $\sigma(S)$, implying that $$\mu_A\big((A\circ \sigma)(S)\big)=\mu_A\big(A(\sigma(S)\big)
=n\mu_A\big(\sigma(S)\big)=n\mu_A(S).$$
Thus, $\mu_A$ is the balanced measure for $A\circ \sigma$, and hence $\mu_A=\mu_{A\circ \sigma}.$ \qed

It was proved by Levin (\cite{lev}, \cite{lev2})  
that  for any rational function  $A$ of degree $n\geq 2$ that is not conjugate to $z^{\pm n}$ there exist at most finitely many rational functions $B$ of any given degree $d\geq 2$ sharing the measure of maximal entropy with $A$. 
Levin's theorem combined with Lemma \ref{fst} and Lemma \ref{sst}  implies the following result. 

\bt \l{ent2} Let $A$ be a rational function of degree $n\geq 2$ that is not conjugate to $z^{\pm n}.$ Then the groups 
$\Aut_{\infty}(A)$ and $\Sigma_{\infty}(A)$ are finite. 
\et
\pr Since $\sigma\in \Aut_{\infty}(A)$ implies that  $A$ and $\sigma^{-1}\circ A\circ \sigma$
share the  measure of maximal entropy by Lemma \ref{fst}, it follows from 
 Levin's theorem that the set of functions \be \l{fset} \sigma^{ -1}\circ A\circ \sigma, \ \ \ \  \sigma\in \Aut_{\infty}(A),\ee is finite. 
On the other hand, the equality  \be \l{sas} \sigma^{-1}\circ A\circ \sigma=\sigma^{\prime -1}\circ A\circ \sigma^{\prime}, \ \ \ \ \ \ \sigma^{\prime}\in \Aut(\C\P^1), \ee  implies that $\sigma^{\prime}\circ \sigma^{-1} \in \Aut(A)$. Thus, the finiteness of set \eqref{fset} implies that there exist $\sigma_1,\sigma_2,\dots,\sigma_l$ such that any $\sigma^{\prime}  \in \Aut_{\infty}(A)$ has the form 
$$\sigma^{\prime}=\widehat\sigma\circ \sigma_k,$$ for some $\widehat\sigma\in \Aut(A)$ and $k,$ $1\leq k \leq l.$   Since $\Aut(A)$ is finite, this implies that $\Aut_{\infty}(A)$ is also finite.

Similarly, it follows from  Lemma \ref{sst} and  Levin's theorem that  
the set of functions $$ A\circ \sigma, \ \ \ \  \sigma\in \Sigma_{\infty}(A),$$ is finite, implying the finiteness of   $\Sigma_{\infty}(A)$ since the equality 
$$A\circ \sigma=A\circ \sigma^{\prime}$$ 
yields that $\sigma^{\prime}\circ \sigma^{-1} \in \Sigma(A)$. \qed

\section{Groups  $\widehat G(A^{\circ k})$ and $\Aut_{\infty}(A)$ }
Let $A$ be a rational function of degree $n\geq 2$. 
We define the set $S(A)$ as the union 
$$S(A)=\bigcup_{i=1}^{\infty} \widehat G(A^{\circ k}),$$ that is, as
the set of M\"obius transformation  $\nu$ such that the equality 
\be \l{egi} \nu \circ A^{\circ k}=A^{\circ k}\circ \mu\ee holds for some  M\"obius transformation $\mu$ and $k\geq 1$. 
The next several results provide a characterization of elements of $S(A)$ and show that $S(A)$ is finite and bounded in terms of $n$, unless $A$ is a quasi-power.   

We start from the following statement.

\bt\l{ker0} 
 Let $A_1,A_2,\dots, A_k$  and  $B_1,B_2,\dots, B_k$, $k\geq 2,$
be rational functions of degree $n\geq 2$ such that 
\be \l{ebs} 
A_1\circ A_2\circ \dots \circ A_k=B_1\circ B_2\circ \dots \circ B_k.
\ee
 Then 
$c(A_1)\subseteq c(B_1\circ B_2)$. 
\et
\pr Let $f$ be a rational function of degree $d$, and $T\subset \C\P^1$ a finite set. 
It is clear that the cardinality of the preimage  $f^{-1}(T)$ satisfies the upper bound 
\be \l{g1} \vert f^{-1}(T)\vert \leq \vert T \vert d.\ee 
To obtain the lower bound, we observe that the Riemann-Hurwitz formula 
$$2d-2=\sum_{z\in \C\P^1}(\mult_zf-1)$$ implies that 
$$\sum_{z\in f^{-1}(T)}(\mult_zf-1)\leq 2d-2.$$  Therefore, 
\be \l{g0} \vert f^{-1}(T)\vert =\sum_{z\in f^{-1}\{T\}} 1 \geq 
\sum_{z\in f^{-1}\{T\}}\mult_zf- 2d+2=(\vert T \vert -2)d+2.\ee

Let us denote by $F$ the rational function defined by any of the parts of equality \eqref{ebs}. 
Assume that $c$ is a critical value of $A_1$ such that $c\notin c(B_1\circ B_2)$. 
Clearly,  
$$\vert F^{-1}\{c\}\vert =\vert (A_2\circ \dots \circ A_k)^{-1}(A_1^{-1}\{c\})\vert.$$ Therefore, since $c\in c(A_1)$ implies that  $\vert A_1^{-1}\{c\}\vert \leq n-1$, it follows from \eqref{g1} that 
\be \l{fr1} \vert F^{-1}\{c\}\vert \leq (n-1)n^{k-1}.\ee On the other hand, 
$$\vert F^{-1}\{c\}\vert =\vert(B_3\circ \dots \circ B_k)^{-1}((B_1\circ B_2)^{-1}\{c\})\vert.$$
Since the condition $c\notin c(B_1\circ B_2)$ is equivalent to the equality 
$\vert (B_1\circ B_2)^{-1}\{c\}\vert =n^2,$  this implies by \eqref{g0} that 
\be \l{fr2} \vert F^{-1}\{c\}\vert \geq (n^2-2)n^{k-2}+2.\ee 
It follows now from \eqref{fr1} and \eqref{fr2} that 
 $$(n^2-2)n^{k-2}+2\leq (n-1)n^{k-1},$$ 
or equivalently that $n^{k-1}+2\leq 2n^{k-2}.$ However, this
leads to a contradiction since $n\geq 2$ implies that $n^{k-1}+2\geq 2n^{k-2}+2.$ Therefore, $c(A_1)\subseteq c(B_1\circ B_2)$.\qed  

\vskip 0.2cm

Theorem \ref{ker0} implies the following statement. 
%, which is essentially the  first statement of Theorem \ref{main}. 

\bt\l{ker} 
 Let $A$ be a rational function of degree $n\geq 2$. Then for every $\nu\in S(A)$ the inclusion 
$\nu\big(c(A)\big)\subseteq c(A^{\circ 2})$ holds. 
\et

\pr 
Let  $\nu$ be an element of $S(A)$. In case  $\nu \in  \widehat G(A)$, the statement of the theorem follows from Lemma \ref{meb}, iii), since $c(A)\subseteq c(A^{\circ 2})$ by the chain rule. Similarly, if $\nu$ belongs to $\widehat G(A^{\circ 2})$, then $\nu\big(c(A^{\circ 2})\big)= c(A^{\circ 2})$, implying that  
$$\nu\big(c(A)\big)\subseteq  \nu\big(c(A^{\circ 2})\big)= c(A^{\circ 2}).$$ Therefore, we may assume that $\nu \in  \widehat G(A^{\circ k})$ for some $k\geq 3.$ Since equality \eqref{egi} has the form  \eqref{ebs}  with 
$$A_1=\nu \circ A, \ \ \ \ \ \ A_2=A_3=\dots= A_k=A, $$
and 
$$B_1=B_2=\dots= B_{k-1}=A, \ \ \ \ \ \  B_k=A\circ \mu,$$  applying Theorem \ref{ker0} 
we conclude that $c(\nu\circ A)\subseteq c(A^{\circ 2})$. Taking into account that for any rational function $A$
 the equality $$c(\nu\circ A)=\nu\big(c(A)\big)$$ holds, this implies that 
$\nu\big(c(A)\big)\subseteq c(A^{\circ 2}).$ \qed

\bt\l{ker2} 
 Let $A$ be a rational function of degree $n\geq 2$. Then the set $S(A)$ is finite and  bounded in terms of $n$, 
unless $A$ is  a quasi-power. Furthermore, the  set $\bigcup_{i=2}^{\infty} \widehat G(A^{\circ k})$ is finite and bounded in terms of $n$, unless $A$ is conjugate to $z^{\pm n}.$ 
\et 
\pr Since any M\"obius transformation is defined by its values at any three points, 
 the condition  $\nu\big(c(A)\big)\subseteq c(A^{\circ 2})$
is satisfied only for finitely many M\"obius transformations whenever $A$ has at least three critical values. Thus, the finiteness of $S(A)$ in case  $A$ is not a quasi-power follows from the first part of Lemma \ref{qua}. 
Moreover, since $\vert c(A)\vert $ and $\vert c(A^{\circ 2})\vert $ are bounded in terms of $n$, the set $S(A)$ is also bounded in terms of $n$.

Further, if $A$ is not conjugate to $z^{\pm n}$, then its second iterate $A^{\circ 2}$ is not a quasi-power by the second part of Lemma \ref{qua}. To prove the finiteness of $\bigcup_{i=2}^{\infty} \widehat G(A^{\circ k})$ in this case, it is enough to show that for every $\nu\in \widehat G(A^{\circ k}),$ $k\geq 2,$ the inclusion   
\be \l{intu} \nu\big(c(A^{\circ 2})\big)\subseteq c(A^{\circ 4})\ee holds, 
 and this can be done by a modification of the proof of Theorem \ref{ker}. 
Indeed, 
equality \eqref{egi} implies the equality 
$$ \nu \circ A^{\circ 2k}=A^{\circ k}\circ \mu \circ A^{\circ k}$$
which can be rewritten for $k\geq 4$ in the form 
 \eqref{ebs} with 
$$A_1=\nu \circ A^{\circ 2}, \ \ \ \ \ \ A_2=A_3=\dots= A_k=A^{\circ 2}, $$
and 
$$B_1=\dots= B_{\frac{k}{2}}=A^{\circ 2}, \ \ \ \   B_{\frac{k}{2}+1}= \mu\circ A^{\circ 2},\ \ \ \ B_{\frac{k}{2}+2}=\dots= B_{k}=A^{\circ 2},$$ 
if $k$ is even, or 
$$B_1=\dots= B_{\frac{k-1}{2}}=A^{\circ 2}, \ \ \ \   B_{\frac{k-1}{2}+1}=A\circ \mu\circ A,\ \ \ \ B_{\frac{k-1}{2}+2}=\dots= B_{k}=A^{\circ 2}, $$ 
if $k$ is odd. Therefore, if $\nu$ belongs to $\widehat G(A^{\circ k})$ for some $k\geq 4,$ then 
applying \linebreak Theorem \ref{ker0}, 
we conclude that \eqref{intu} holds. 
On the other hand, 
if $\nu$ belongs to $\widehat G(A^{\circ 2})$, then $\nu\big(c(A^{\circ 2})\big)= c(A^{\circ 2})$, by Lemma \ref{meb}, iii), implying \eqref{intu} by the chain rule. 
Similarly, if $\nu$ belongs to $\widehat G(A^{\circ 3})$, then $\nu\big(c(A^{\circ 3})\big)= c(A^{\circ 3})$, implying that  
$$\nu\big(c(A^{\circ 2})\big)\subseteq  \nu\big(c(A^{\circ 3})\big)= c(A^{\circ 3})\subseteq c(A^{\circ 4}).\eqno{\Box}$$

 Theorem \ref{ker2} implies the following result.

\bt \l{sig33} Let $A$ be a rational function of degree $n\geq 2$.
Then   the orders of the groups   $\widehat G(A^{\circ k})$, $k\geq 1,$ 
are finite and uniformly bounded in terms of  $n$ only, unless $A$ is  a quasi-power. 
Furthermore, the orders of the groups   $\widehat G(A^{\circ k})$, $k\geq 2,$ 
are finite and uniformly bounded in terms of  $n$ only, unless $A$ is conjugate to $z^{\pm n}.$ 
\et
\pr 
%Since every group $\widehat G(A^{\circ k})$, $k\geq 1,$ is contained in $S(A)$, while 
%every group $\widehat G(A^{\circ k})$, $k\geq 2,$ is contained in  $S(A)\setminus \widehat G(A)$ 
The theorem is a direct corollary of  Theorem \ref{ker2}. \qed  

Finally, Theorem \ref{ker} and Theorem \ref{ker2} imply Theorem \ref{main} from the introduction.

\vskip 0.2cm

\noindent{\it Proof of Theorem \ref{main}.} 
The boundedness of the set 
$\bigcup_{i=2}^{\infty} \Aut(A^{\circ k})$ in terms of $n$ for $A$ that is not conjugate to $z^n$ follows from Theorem \ref{ker2}. On the other hand, $\Aut(A)$ is finite and bounded in terms of $n$ by Corollary \ref{yc2}. This proves the first part of the theorem. 
Finally, since the set $S(A)$ contains the group  $\Aut_{\infty}(A)$, the second part of the theorem follows from Theorem \ref{ker} (the assumption that $A$ is not conjugate to $z^n$ is actually redundant).    \qed

\section{Groups $\Sigma_{\infty}(A)$ and $G(A^{\circ k})$} 

Let $A$ and $B$ be  rational functions of degree at least two. We 
recall that the function $B$ is said to be {\it semiconjugate} to the function $A$
if there exists a non-constant rational function $X$
such that the equality 
\be \l{i1} A\circ X=X\circ B\ee
holds. Usually, we will write this condition in the form of a commuting diagram 
\be 
\begin{CD} 
\C\P^1 @> B>>\C\P^1 \\ 
@V X  VV @VV X  V\\ 
 \C\P^1 @> A>> \C\P^1.
\end{CD} 
\ee
The simplest examples of semiconjugate rational functions are provided by equivalent rational functions  defined in the introduction. Indeed, it follows from equalities 
\eqref{ey} that the diagrams 
\be 
\begin{CD} 
\C\P^1 @>  A>>\C\P^1 \\ 
@V {V}  VV @VV {V}  V\\ 
 \C\P^1 @> \tt A>> \C\P^1
\end{CD} \ \ \ \ \ \ \  \ \ \ \ 
\begin{CD} 
\C\P^1 @> \tt A>>\C\P^1 \\ 
@V U  VV @VV U  V\\ 
 \C\P^1 @> A>> \C\P^1
\end{CD} 
\ee
commutes, implying inductively that if $A$ is equivalent to $B$, then $A$ is semiconjugate to $B$, and  $B$ is semiconjugate to $A$. 

A comprehensive description of semiconjugate  rational functions was obtained in the  papers  \cite{semi}, \cite{rec}, \cite{lattes}.
In particular, it was shown in \cite{semi} that solutions $A,X,B$ 
of \eqref {i1} satisfying  
$\C(X,B)=\C(z),$  called {\it primitive}, can be described in terms of group actions on $\C\P^1$ or $\C$, implying strong restrictions on a possible form of $A$, $B$ and $X$. 
On the other hand, an arbitrary solution of equation \eqref{i1} can be reduced to a primitive one by a sequence of elementary transformations as follows.  By the L\"uroth theorem, the field 
$\C(X,B)$ is generated by some rational function $W$.  
Therefore, if  $\C(X,B)\neq \C(z)$, then 
 there exists a rational function $W$ of degree greater than one such that
\be \l{of} B=\tt B\circ W, \ \ \ X=\tt X\circ W\ee 
for some rational functions $\tt X$ and $\tt B$ satisfying $\C(\tt X,\tt B)=\C(z)$. Moreover, it is easy to see  that the diagram
\be  
\begin{CD} 
\C\P^1 @> B>>\C\P^1 \\ 
@V W  VV @VV W   V\\ 
\C\P^1 @> W\circ \tt B >>\C\P^1 \\ 
@V \tt X  VV @VV \tt X  V\\ 
 \C\P^1 @> A>> \C\P^1
\end{CD} 
\ee
 commutes. Thus, the triple $A, \tt X,W\circ \tt B$ is another solution of \e{i1}. This new solution is not necessarily primitive, however 
$\deg \tt X<\deg X$. Therefore, continuing in this way, after a finite number of similar transformations we will arrive to a primitive solution. In more detail, the above argument shows  that for any rational functions $A,X,B$ satisfying \eqref{i1}  there exist rational functions $X_0$, $B_0$,   $U$ such that $X=X_0\circ U,$ the   diagram 
\be \l{diag} 
\begin{CD} 
\C\P^1 @> B>>\C\P^1 \\ 
@V U  VV @VV U   V\\ 
 \C\P^1 @> B_0>> \C\P^1
 \\ 
@V {X_0}  VV @VV {X_0}   V\\ 
 \C\P^1 @> A>> \C\P^1 
\end{CD} 
\ee
commutes,  the triple 
$A,X_0, B_0$ is a primitive solution of \eqref{i1}, and $B_0\sim B$.

The following theorem is essentially the second part of Theorem \ref{sug} from the introduction but without the assumption that $A$ is not conjugate to $z^n$, which is redundant in this case.

\bt \l{sig} Let $A$ be a rational function of degree $n\geq 2$. Then 
for every $\sigma\in\Sigma_{\infty}(A)$ the relation $A\circ \sigma\sim A$ holds. 
\et
\pr Let $\sigma$ be an element of $
\Sigma_{\infty}(A)$. Then \be \l{kotik} A^{\circ k}=A^{\circ k}\circ \sigma\ee for some $k\geq 1.$ 
Writing this equality 
as the semiconjugacy 
\be \l{krol} 
\begin{CD}
\C\P^1 @> A\circ \sigma >> \C\P^1\\
@VV A^{\circ (k-1)} V @VV A^{\circ (k-1)} V\\ 
\C\P^1 @>A >>\C\P^1\,, 
\end{CD}
\ee
we see that to prove  the theorem it is enough to show that in  diagram 
\eqref{diag}, corresponding to the solution 
$$A=A,\ \ \ \ X=A^{\circ (k-1)},\ \ \ \ B=A\circ \sigma$$ of \eqref{i1},  
the function  $X_0$ has degree one. The proof of the last statement is similar to the proof of Theorem 2.3 in \cite{rev} and follows from the following two facts. First, for any   primitive solution $A,X,B$ of \eqref{i1}, the solution $A^{\circ l},X,B^{\circ l}$, $l\geq 1$, is also primitive (see \cite{rev}, Lemma 2.5). Second, a solution $A,X,B$ of \eqref{i1} is primitive if and only if the algebraic curve $$A(x)-X(y)=0$$ is irreducible (see \cite{rev}, Lemma 2.4).
Using these facts we see that  
 the triple $A^{\circ (k-1)},X_0,B_0^{\circ (k-1)}$  is a primitive solution of \eqref{i1}, and the algebraic curve \be \l{cucu} A^{\circ (k-1)}(x)-X_0(y)=0\ee is irreducible. However, the equality $$A^{\circ (k-1)}=X_0\circ U,$$ implies that the curve 
$$U(x)-y=0$$ is a component of \eqref{cucu}. Moreover, if $\deg X_0>1$, then this component is proper. Therefore, $\deg X_0=1$. \qed

The following result proves the first part of Theorem \ref{sug} and thus finishes the proof of this theorem. 

\bt \l{sig2}  Let $A$ be a rational function of degree $n\geq 2$ that  is not conjugate to $z^{\pm n}.$  Then the order of the group  $\Sigma_{\infty}(A)$ is finite and  bounded in terms of $n$.
\et
\pr
Let us observe first that it is enough to prove the theorem under the assumption that $A$ is not a quasi-power.  Indeed, if $A$ is a quasi-power but is not conjugate to $z^{\pm n},$ then $A^{\circ 2}$ is not a quasi-power by Lemma \ref{qua}. Therefore, if the theorem is true for functions that are not quasi-powers, then for any $A$ that is not conjugate to  $z^{\pm n}$, the group 
$\Sigma_{\infty}(A^{\circ 2})$ is finite and  bounded in terms of $n$, implying by \eqref{ospes} that the same is true for the group $\Sigma_{\infty}(A)$.

Assume now that $A$ is not a quasi-power. Then $G(A)$ is finite by Theorem \ref{prim}. Let us recall that in view of equality \eqref{lus} the equivalence class $[A]$ is a union of conjugacy classes. 
Denoting the number of these conjugacy classes by $N_A$, let us show that  if $N_A$ 
is finite, then \be \l{ine} \vert \Sigma_{\infty}(A)\vert \leq \vert G(A)\vert N_A. \ee
 By Theorem \ref{sig},  for any $\sigma\in \Sigma_{\infty}(A)$ the function 
$A\circ \sigma$ belongs to one of $N_A $ conjugacy classes in the equivalence class $[A]$. 
Furthermore, if $A\circ \sigma_0$ and  $A\circ \sigma$ 
belong to the same conjugacy class, then 
$$A\circ \sigma=\alpha\circ A\circ \sigma_0\circ \alpha^{-1}$$ for some $\alpha \in \Aut(\C\P^1),$ 
implying that  
$$A\circ \sigma\circ \alpha\circ \sigma_0^{-1}=\alpha\circ A.$$
This is possible only if $\alpha$ belongs to the group $\widehat G(A)$, and, in addition,  
$\sigma\circ \alpha\circ \sigma_0^{-1}$ belongs to the preimage of $\alpha$ under homomorphism \eqref{homic}.  Therefore, 
for any fixed $\sigma_0$, there could be at most 
$\vert \widehat G(A)\vert$ such $\alpha$, and for each  $\alpha$ there could be at most $\vert \Ker \gamma_A\vert$ elements $\sigma\in \Sigma_{\infty}(A)$ such that  
$$\gamma_A(\sigma\circ \alpha\circ \sigma_0^{-1})=\alpha.$$ Thus, \eqref{ine} follows from \eqref{burn}.

It was proved in \cite{rec} that  $N_A$ is infinite  if and only if $A$ is a flexible Latt\`es map. However, the proof given in \cite{rec} uses the theorem of McMullen (\cite{Mc}) about isospectral rational functions, which is not effective. Therefore,  the result of  \cite{rec} does not imply that $N_A$ is bounded in terms of $n$. Nevertheless, we can  use the main result of  \cite{fin}, which yields in particular that for a given rational function $B$ of degree $n\geq 2$ the number of conjugacy classes of rational functions $A$ such that \eqref{i1} holds for some rational function $X$ is finite and bounded in terms of $n$, unless $B$  is {\it special}, that is, unless $B$ is either  a Latt\`es map or it is
conjugate to $z^{\pm n}$ or $\pm T_n.$
Since $A\sim A'$ implies that $A$ is semiconjugate to $A'$, this  implies that for non-special $A$ the number $N_A$ is bounded in terms of $n$. Moreover, it is easy to see that the same is  true also for $A$ 
conjugate to $z^{\pm n}$ or $\pm T_n,$ since any decomposition of $z^n$ has the form 
$$z^n=(z^d\circ \mu)\circ (\mu^{-1}\circ z^{n/d}),$$ where $\mu\in \Aut(\C\P^1)$ and $d\vert n,$ while 
any decomposition of $T_n$ has the form 
$$T_n=(T_d\circ \mu)\circ (\mu^{-1}\circ T_{n/d}),$$ where $\mu\in \Aut(\C\P^1)$ and $d\vert n$. 

The above shows that to finish the proof of Theorem \ref{sig2} we only must prove that the group $\Sigma_{\infty}(A)$ is finite and bounded in terms of $n$ if $A$ is a Latt\`es map.
% or is conjugate to $\pm T_n$.  
%It is easy to see using the explicit formula 
%\be \l{cheb} T_n=\frac{n}{2}\sum_{k=0}^{[n/2]}(-1)^k\frac{(n-k-1)!}{k!(n-2k)!}(2x)^{n-2k}\ee 
% that the group $\Sigma(\pm T_n)$ is either trivial or equal to $C_2,$ depending on the parity of $n$. Therefore, since $T_n^{\circ k}=T_{n^{\circ k}},$ $k\geq 1,$  the order of $\Sigma_{\infty}(\pm T_n)$ is at most two. 
To prove the last statement, we recall that if $A$ is a Latt\`es map, then there exists an orbifold  $\f O=(\C\P^1,\nu)$  of zero Euler characteristic  such that $A:\, \f O\rightarrow \f O$ is a covering map between orbifold (see \cite{mil}, \cite{lattes} for more detail).
Since this implies that $A^{\circ k}:\, \f O\rightarrow \f O$, $k\geq 1$, also is a covering map (see \cite{semi}, Corollary 4.1), it follows from equality \eqref{kotik}  that  $\sigma:\, \f O\rightarrow \f O$ is a covering map (see \cite{semi}, Corollary 4.2 and Lemma 4.1). As $\sigma$ is of degree one, the last condition simply means that $\sigma$ permute points of the support of $\f O$. Since the support of an orbifold $\f O=(\C\P^1,\nu)$ of zero Euler characteristic contains either three or four points, this implies that  
$\Sigma_{\infty}(A)$ is finite and uniformly bounded for any Latt\`es map $A$. 
%This finishes the proof.
\qed 

\vskip 0.2cm

\noindent{\it Proof of Theorem \ref{sig4}.} If $\sigma\in \Sigma_{\infty}(A)$, then \be \l{kro} A\circ \sigma\sim A,\ee  by Theorem \ref{sig}. On the other hand, since for any indecomposable function $A$ 
 the number $N_A$ obviously is equal to one,   condition  
\eqref{kro} is equivalent to the condition that
\be \l{kross} A\circ \sigma=\beta \circ A \circ \beta^{-1}\ee for some $\beta \in \Aut(\C\P^1)$. 
Clearly, equality 
\eqref{kross} implies that $\beta$ belongs to $\widehat G(A)$. Therefore, if  $\widehat G(A)$ is trivial, then \eqref{kro} is satisfied only if $A\circ \sigma= A$, that is, only if $\sigma$ belongs to $\Sigma(A)$. Thus, 
$\Sigma(A)=\Sigma_{\infty}(A)$, whenever  $\widehat G(A)$ is trivial.  

Furthermore,  it follows from equality \eqref{kross}  that  $\sigma\circ \beta$ belongs to the preimage of $\beta$ under homomorphism \eqref{homic}. On the other hand,  if $G(A)=\Aut(A)$, this preimage consists of $\beta$ only. Therefore, in this case   $\sigma\circ \beta=\beta,$ implying that $\sigma$ is the identity map. Thus,  
the group $\Sigma_{\infty}(A)$ is  trivial, whenever $G(A)=\Aut(A)$.  \qed

The following  theorem implies Theorem \ref{sig0}  from the introduction. 

\bt \l{sig3} Let $A$ be a rational function of degree $n\geq 2$.
Then   the orders of the groups   $ G(A^{\circ k})$, $k\geq 1,$ 
are finite and uniformly bounded in terms of  $n$ only, unless $A$ is  a quasi-power. 
Furthermore, the orders of the groups   $ G(A^{\circ k})$, $k\geq 2,$ 
are finite and uniformly bounded in terms of  $n$ only, unless $A$ is conjugate to $z^{\pm n}.$ 
\et
\pr If $A$ is not a quasi-power, then by Theorem \ref{sig33} and Theorem \ref{sig2}   
the orders of the groups  $\widehat G(A^{\circ k})$, $k\geq 1,$ and  $\Sigma(A^{\circ k})$, 
$k\geq 1,$  are finite and uniformly bounded  in terms of  $n$ only. 
Therefore, by \eqref{burn},  the orders of the groups  $ G(A^{\circ k})$, $k\geq 1,$
also are finite and uniformly bounded.   
Similarly, the groups    $ G(A^{\circ k})$, $k\geq 2,$ 
are finite and uniformly bounded  in terms of  $n$ only, unless $A$ is conjugate to $z^{\pm n}.$ \qed

\bc\l{ya}  Let $A$ be a rational function of degree $n\geq 2$. Then the sequence $ G(A^{\circ k})$, $k\geq 1,$  
contains only finitely many non-isomorphic groups. 
\ec 
\pr 
For $A$ not  conjugate to $z^{\pm n}$,  the corollary follows from Theorem \ref{sig3}  
since there exist only finitely many groups of any given order. Moreover,  actually  
the groups $ G(A^{\circ k})$, $k\geq 2,$ belong to the list $A_4,$ $S_4,$ $A_5,$ $C_l$, $D_{2l}$,  
 by Theorem \ref{prim}. 
On the other hand, if  $A$ is conjugate to $z^{\pm n}$, then all the groups $G(A^{\circ k})$, $k\geq 1,$ consist of the transformations $z\rightarrow c z^{\pm 1}$, $c\in \C\setminus\{0\}$.
 \qed

We finish this section by two examples of calculation of the group $\Sigma_{\infty}(A)$.

\begin{example}
Let us consider the function $$A=x+\frac{27}{x^{3}}.$$ A calculation shows that, in addition to the critical value $\infty$, this function has critical values $\pm 4$ and $ \pm 4i,$ and $$A\pm 4={\frac { \left( {x}^{2}\mp 2\,x+3 \right)  \left( x\pm 3 \right) ^{2}}{{x}^{
3}}},$$ $$ A\pm 4i={\frac { \left( {x}^{2}\mp 2\,ix-3 \right)  \left( \pm x+3\,i \right) ^{2}}{
{x}^{3}}}.$$ Since the above equalities imply that $\mult_0A=3$, while at any other point of $\C\P^1$ the multiplicity of $A$ is at most two, it follows from Theorem \ref{cyc} that $G(A)$ is a cyclic group, whose generator has zero as a fixed point. Moreover, since $G(A)$ obviously contains the transformation $\sigma=-z$, the second fixed point of this generator must be infinity. This  
implies easily that $G(A)$ is a cyclic group of order two, and  $G(A)= \Aut(A)$. Finally,  since $\mult_0A=3$, 
it follows from the chain rule that the equality $A=A_1\circ A_2$, where $A_1$ and $A_2$ are rational function of degree two is impossible. Therefore, $A$ is indecomposable, and hence the group $\Sigma_{\infty}(A)$ is trivial by Theorem \ref{sig4}.
\end{example}

\begin{example} \l{notdef0} 
Let us consider  the function 
$$A=\frac{z^2-1}{z^2+1}.$$ Since $A$ is a quasi-power, $\Sigma(A)$ is a cyclic group of order two, generated by the transformation $z\rightarrow -z.$ A calculation shows that the second iterate  
$$A^{\circ 2}=-\frac{2z^2}{z^4+1}$$ is the function $B$ from Example \ref{dadef}. Thus,  
 $\Sigma(A^{\circ 2})$ is the dihedral group $D_4$, generated by the transformation $z\rightarrow -z$ and $z\rightarrow 1/z$. In particular, $\Sigma(A^{\circ 2})$  is larger than $\Sigma(A)$.  Moreover, since 
$$A^{\circ 3}=-{\frac { \left( {z}^{4}-1 \right) ^{2}}
{{z}^{8}+6\,{z}^{4}+1}}
,$$  we see that $\Sigma(A^{\circ 3})$  contains the dihedral group $D_8$, generated by the transformation $\mu_1=iz$ and $\mu_2= 1/z,$ and hence  $\Sigma(A^{\circ 3})$  is larger than  $\Sigma(A^{\circ 2})$.

Let us show that $$\Sigma_{\infty}(A)=\Sigma(A^{\circ 3})=D_8.$$ As in Example \ref{dadef}, we see that if $\Sigma_{\infty}(A)$ is larger than $D_8$, then either  $\Sigma_{\infty}(A)=S_4$, or  $\Sigma_{\infty}(A)$ is a dihedral group containing an element $\sigma$ of order $l>4$
such that $\mu_1$ is an iterate of $\sigma.$ 
The first case is  impossible, for   
 otherwise Theorem \ref{aut} implies that for $k$ satisfying $\Sigma_{\infty}(A)=\Sigma(A^{\circ k})$ the number 
$\deg A^{\circ k}=2^k$ is divisible by  $\vert S_4\vert =24.$ On the other hand,  in the second case, the fixed points of $\sigma$ are zero and infinity. 
%Therefore, 
%by Theorem \ref{sig}, taking into account that $A$ is indecomposable, 
Since  $A$ is indecomposable, it follows from Theorem \ref{sig} that 
to exclude the second case it  is enough to show that  if  $\sigma=cz,$ $c\in \C\setminus\{0\},$ satisfies 
 \be \l{koro} A\circ \sigma=\beta \circ A \circ \beta^{-1}, \ \ \  \beta\in \Aut(\C\P^1),\ee
 then $\sigma$ is an iterate of $\mu_1$.
Since critical points of the function on the left side of \eqref{koro} coincide with 
 critical points of the function on the right side, the M\"obius transformation  $\beta$ 
 necessarily has the form  $\beta=d z^{\pm 1}$, $d\in \C\setminus\{0\}.$
Thus, equation \eqref{koro} reduces to the equations $$\frac{c^2z^2-1}{c^2z^2+1}=\frac{1}{d}\frac{d^2z^2-1}{d^2z^2+1} ,$$ and 
$$\frac{c^2z^2-1}{c^2z^2+1}=
{\frac {d \left( {d}^{2}+{z}^{2} \right) }{{d}^{2}-{z}^{2}}}.$$
One can check that solutions of the first equation are $d=1$ and $c=\pm 1$, while solutions of the second are 
$d=-1$ and $c=\pm i.$ This proves the necessary statement. 
Notice that instead of Theorem \ref{sig} it is also possible to use Theorem \ref{seq0} (see the next section).

\end{example}

\section{Groups $G(A,z_0,z_1)$}
Following \cite{rez}, we say that 
a formal power series $f(z)=\sum_{i=1}^{\infty}a_iz^i$ having zero as a fixed point is {\it homozygous} $\mod l$ if the inequalities $a_i\neq 0$ and $a_j\neq 0$ imply the equality $i\equiv j (\mod l).$ If $f$ is not homozygous  $\mod l$, it is called {\it hybrid} $\mod l.$ 
Obviously,  the condition that $f$ is homozygous  $\mod l$ 
is equivalent to the condition that  $f=z^rg(z^l)$ for some formal power series $g=\sum_{i=0}^{\infty}b_iz^i$  and integer $r,$ \linebreak $1\leq r\leq l$. In particular,   if $f$ is  homozygous $\mod l$, then  any iterate of $f$ is
homozygous $\mod l$.  The inverse is not true. However, the following statement  proved by 
 Reznick (\cite{rez}) holds:  if a formal power series $f(z)=\sum_{i=1}^{\infty}a_iz^i$ is hybrid $\mod l$ and $f^{\circ k}$ is homozygous $\mod l$, then $f^{\circ ks}(z)=z$ for some integer $s\geq 1$. Our proof of 
 Theorem \ref{seq0} relies on this result.

\vskip 0.2cm

\noindent {\it Proof of Theorem \ref{seq0}.} Without loss of generality, we can assume that $z_0=0$ and $z_1=\infty$. 
Let  $f_A$ be the Taylor series of the function $A$ at zero. Arguing as in the proof of Theorem \ref{prim}, we see that  every element of $G(A,0,\infty)$ has the form $z\rightarrow \v z$, where $\v$ is a root of unity, and  $G(A,0,\infty)$ 
is a finite cyclic group, whose order is equal to the maximum number $n$ such that 
$f_A$ is homozygous $\mod n$. Since  $f_{A^{\circ k}}=f_A^{\circ k},$  this implies that 
$$ G(A,0,\infty)\subseteq  G(A^{\circ k},0,\infty), \ \ \ \ \  k\geq 1.$$ 
Moreover, if $G(A^{\circ k},0,\infty)$ is strictly larger than $G(A,0,\infty)$ for some $k>1$, then there exists $n_0$ such that $f_A$ is hybrid $\mod n_0$ but $f_A^{\circ k}$ is homozygous $\mod n_0$. Therefore,  by the Reznick theorem,   
the equality $f_A^{\circ ks}=z$ holds for some $s\geq 1$. 
However, in this case by the analytical continuation   $A^{\circ ks}=z$ for all $z\in \C\P^1$,  in contradiction with $n\geq 2$. Thus, the groups $G(A^{\circ k},0,\infty)$, $k\geq 1,$ are  equal. \qed 

Notice that  the groups $G(A^{\circ k},z_0,z_1)$, $k\geq 1$, are equal even if $A$ is conjugate to $z^{\pm n}$. 
Indeed, for $A=z^{\pm n}$ these groups are trivial, unless $\{z_0,z_1\}=\{0,\infty\}$, while in the last case all these groups consist of the transformations $z\rightarrow c z^{\pm 1}$, $c\in \C\setminus\{0\}$.

Let us emphasize that since iterates $A^{\circ k }$, $k>1$, have in general more fixed points than $A$, it may happen that 
 $G(A^{\circ k},z_0,z_1)$, $k>1$, is non-trivial, while  $G(A,z_0,z_1)$ is  not defined, so that  
the equality $G(A^{\circ k},z_0,z_1)=G(A,z_0,z_1)$ does not make sense.
For example, for  
the function $$A=\frac{z^2-1}{z^2+1}$$ from Example \ref{notdef0}  zero is not a fixed point for $A$, and hence the group $G(A,0,\infty)$ is not defined. However,  zero is a fixed point for 
$$A^{\circ 2}=-\frac{2z^2}{z^4+1},$$ and the group $G(A^{\circ 2},0,\infty)$ is a cyclic group of order four. 
Let us remark that Theorem \ref{seq0} gives  
 another proof of the fact that $\Sigma_{\infty}(A)$  cannot contain an element $\sigma=c z$, $c\in \C\setminus\{0\}$, of order $l>4$.  Indeed, such $\sigma$  must belong to the group  $G(A^{\circ k},0,\infty)$ for some $k\geq 1$, and hence to the group  $G(A^{\circ 2k},0,\infty)$. However, $G(A^{\circ 2k},0,\infty)$ is equal to 
$G(A^{\circ 2},0,\infty)=C_4$ by Theorem \ref{seq0} applied to $A^{\circ 2}.$

Under certain conditions, Theorem \ref{seq0} permits  to estimate the order of the groups $\Aut_{\infty}(A)$ and $\Sigma_{\infty}(A)$ and even to describe these groups explicitly.

\bt \l{qes0} 
Let $A$ be a rational function of degree $n\geq 2$ that is not conjugate to $z^{\pm n}.$ 
Assume that  for some $k\geq 1$ the group $\Aut(A^{\circ k})$ contains an element $\sigma$ of order at least six with fixed points $z_0$ and $z_1$ such that $z_0$ is a fixed point of $A^{\circ k}$.  Then the inequality $\vert \Aut_{\infty}(A)\vert \leq 2\vert G(A^{\circ k},z_0,z_1)\vert$ holds.  Similarly, if $\sigma$ as above is contained in $\Sigma(A^{\circ k})$, then  
$\vert \Sigma_{\infty}(A)\vert  \leq 2\vert G(A^{\circ k},z_0,z_1)\vert$.
\et
\pr Since the maximal order of a cyclic subgroup in the groups $A_4,$ $S_4,$ $A_5$ is five, 
it follows from Corollary \ref{yc2} that if $\Aut(A^{\circ k})$  contains an element $\sigma$ of order $r>5$,  then either $\Aut_{\infty}(A)=C_{s}$ or $\Aut_{\infty}(A)=D_{2s}$, where $r\vert s$. 
Moreover, if  $\sigma_{\infty}$ is an element of order $s$ in $\Aut_{\infty}(A)$, then 
$\sigma$ is an iterate of $\sigma_{\infty}$. In particular, 
fixed points of $\sigma_{\infty}$ coincide with fixed points of $\sigma.$

To prove the  theorem, we only must show that the inequality \be \l{ih} s> \vert G(A^{\circ k},z_0,z_1)\vert\ee is impossible. Assume the inverse.  Since $\sigma_{\infty}$ belongs to $\Aut(A^{\circ k'})$ for  some $k'\geq 1$, it belongs to $\Aut(A^{\circ kk'})$ and   $
G(A^{\circ kk'},z_0,z_1).$ Therefore, if \eqref{ih} holds,  then the group $G(A^{\circ kk'},z_0,z_1)$ 
contains an element of order greater than $\vert G(A^{\circ k},z_0,z_1)\vert $,  in contradiction with the equality  $$G(A^{\circ kk'},z_0,z_1)=G(A^{\circ k},z_0,z_1),$$ provided by Theorem \ref{seq0} applied to $G(A^{\circ k}).$ 
The proof of the inequality for $ \vert\Sigma_{\infty}(A)\vert $ is similar. \qed

\begin{example} Let us consider the function 
$$A=z\frac{z^6-2}{2z^6-1}.$$ It is easy to see that $\Aut(A)$ contains the dihedral group $D_{12}$  generated by the transformations $$z\rightarrow e^\frac{2\pi i}{6}z, \ \ \ \ z\rightarrow 1/z.$$ 
Since zero  is a fixed point of $A$ and $G(A,0,\infty)=C_6$, it follows from Theorem \ref{qes0} that $$ \Aut_{\infty}(A)=\Aut(A)=D_{12}.$$ 
\end{example}

Although the group $\Aut(A^{\circ k})$ does not necessarily contain an element that belongs to $G(A^{\circ k},z_0,z_1),$ it always contains an element that belongs to $G(A^{\circ 2k},z_0,z_1).$ 
More generally, the  following statement holds. 

\bl \l{kott} Let $A$ be a rational function of degree $n\geq 2$, and  
 $\sigma\notin \Sigma(A^{\circ k})$ a M\"obius transformation such  that the   equality 
\be \l{au1} A^{\circ k}\circ \sigma=\sigma^{\circ l}\circ A^{\circ k},\ee holds for some $l\geq 1.$   Then at least one of  the fixed points $z_0,$ $z_1$ of $\sigma$ is a fixed point of $A^{\circ 2k},$ and if $z_0$ is such a point, then  
$\sigma \in  G(A^{\circ 2k},z_0,z_1).$ 
\el 
\pr Clearly, equality \eqref{au1} implies the equalities 
$$\sigma^{\circ l}(A^{\circ k}(z_0))=A^{\circ k}(z_0), \ \ \ 
\sigma^{\circ l}(A^{\circ k}(z_1))=A^{\circ k}(z_1).$$  However, since $\sigma^{\circ l}$ is not the identity map, it has only two  fixed points $z_0,z_1$. Therefore, 
$A^{\circ k}\{z_0,z_1\}\subseteq \{z_0,z_1\},$ implying that at least one of the points $z_0,z_1$ is a fixed point of $A^{\circ 2k}.$ Finally, if $z_0$ is such a point, then  
$\sigma \in  G(A^{\circ 2k},z_0,z_1)$. \qed

Combining Theorem \ref{qes0}  with Lemma \ref{kott} we obtain the following result. 
\bt \l{qes} 
Let $A$ be a rational function of degree $n\geq 2$ that is not conjugate to $z^{\pm n}.$ 
Assume that  for some $k\geq 1$ the group $\Aut(A^{\circ k})$ contains an element $\sigma$ of order at least six with fixed points $z_0,z_1$.
Then $\vert \Aut_{\infty}(A)\vert \leq 2\vert G(A^{\circ 2k},z_0,z_1)\vert$, where $z_0$ is a fixed point of $\sigma$ that is also a fixed point of $A^{\circ 2k}$. \qed 
\et

\end{document}